\documentclass[12pt]{article}  
\def\sq{\hbox {\rlap{$\sqcap$}$\sqcup$}}
\def\sq{\hbox {\rlap{$\sqcap$}$\sqcup$}}
\def\R{ {\rm R \kern -.31cm I \kern .15cm}}
\def\C{ {\rm C \kern -.15cm \vrule width.5pt \kern .12cm}}
\def\Z{ {\rm Z \kern -.27cm \angle \kern .02cm}}
\def\N{ {\rm N \kern -.26cm \vrule width.4pt \kern .10cm}}
\def\1{{\rm 1\mskip-4.5mu l} }
\def\lsim{\raise0.3ex\hbox{$<$\kern-0.75em\raise-1.1ex\hbox{$\sim$}}}
\def\gsim{\raise0.3ex\hbox{$>$\kern-0.75em\raise-1.1ex\hbox{$\sim$}}}
\def\noi{\noindent}

\def\beq{\begin{equation}}   \def\eeq{\end{equation}}
\def\bea{\begin{eqnarray}}  \def\eea{\end{eqnarray}}

\def\noi{\noindent}
\def\beeq{\begin{eqnarray}} \def\eeeq{\end{eqnarray}}
\newcommand\mysection{\setcounter{equation}{0}\section}

\newcounter{hran}

\begin{document} 
\centerline{\Large\bf Uniqueness at infinity in time for the} 
 \vskip 3 truemm 
 \centerline{\Large\bf Maxwell-Schr\"odinger system with arbitrarily} 
 \vskip 3 truemm
 \centerline{\Large\bf  large asymptotic data
 } \vskip 0.8 truecm

\centerline{\bf J. Ginibre}
\centerline{Laboratoire de Physique Th\'eorique\footnote{Unit\'e Mixte de
Recherche (CNRS) UMR 8627}}  \centerline{Universit\'e de Paris XI, B\^atiment
210, F-91405 Orsay Cedex, France}
\vskip 3 truemm

\centerline{\bf G. Velo}
\centerline{Dipartimento di Fisica, Universit\`a di Bologna}  \centerline{and INFN, Sezione di
Bologna, Italy}

\vskip 1 truecm

\begin{abstract}
We prove the uniqueness of solutions of the Maxwell-Schr\"odinger system with given asymptotic behaviour at infinity in time. The assumptions include suitable restrictions on the growth of solutions for large time and on the accuracy of their asymptotics, but no restriction on their size. The result applies to the solutions with prescribed asymptotics constructed in a previous paper.
\end{abstract}

\vskip 1.5 truecm
\noi MS Classification :    Primary 35P25. Secondary 35B40, 35Q40.\par \vskip 2 truemm

\noi Key words : Long range scattering, uniqueness of solutions, Maxwell-Schr\"odinger system. \par 
\vskip 1 truecm

\noindent LPT Orsay 07-38\par
\noindent June 2007\par \vskip 3 truemm

\newpage
\pagestyle{plain}
\baselineskip 18pt
\mysection{Introduction}
\hspace*{\parindent} This paper is a sequel to a previous paper \cite{5r}, hereafter referred to as II, where we studied the theory of scattering for the Maxwell-Schr\"odinger system (MS) in $3 + 1$ dimensional
space time. That system describes the evolution of a charged nonrelativistic quantum
mechanical particle interacting with the (classical) electromagnetic field it generates. It can be written as follows~:
\beq
\label{1.1e}
\left \{ \begin{array}{l} i\partial_t u = -(1/2) \Delta_A u + 
A_{e} u \\ \\ \sq A_{e} - \partial_t \left ( \partial_t A_{e} + \nabla \cdot A\right )  = |u|^2 \\ \\ \sq A + \nabla \left ( \partial_t A_{e} + \nabla \cdot A\right )  = {\rm Im}\ \overline{u} \nabla_A u\end{array} \right . \eeq

\noi where $u$ and $(A, A_e)$ are respectively a complex valued function and an ${I\hskip-1truemm R}^{3+1}$ valued function defined in space time ${I\hskip-1truemm R}^{3+1}$, $\nabla_A = \nabla - iA$ and $\Delta_A = \nabla_A^2$ are the covariant gradient and covariant Laplacian respectively, and $\sq = \partial_t^2 - \Delta$ is the d'Alembertian. An important property of that system is its gauge invariance, namely the invariance under the transformation
$$\left ( u, A, A_e \right ) \rightarrow \left (u \exp (- i \theta), A - \nabla \theta, A_e + \partial_t \theta \right ) \ ,$$

\noi where $\theta$ is an arbitrary real function defined in ${I\hskip-1truemm R}^{3+1}$. As a consequence of that invariance, the system (\ref{1.1e}) is underdetermined as an evolution system and has to be supplemented by an additional equation, called a gauge condition. In this paper,  we shall use exclusively the Coulomb gauge condition, namely $\nabla \cdot A = 0$. Under that condition, the equation for $A_e$ can be solved by
\beq
\label{1.2e}
A_e = - \Delta^{-1} |u|^2 = (4 \pi |x|)^{-1} \star |u|^2 \equiv g(u)
\eeq

\noi where $\star$ denotes the convolution in ${I\hskip-1truemm R}^3$. Substituting (\ref{1.2e}) and the gauge condition into (\ref{1.1e}) yields the formally equivalent system
\bea
\label{1.3e}
&&i \partial_t u = - (1/2) \Delta_A u + g(u) u \\
&&\sq A = P \ {\rm Im} \ \overline{u} \nabla_A u
\label{1.4e}
\eea

\noi where $P = \1 - \nabla \Delta^{-1} \nabla$ is the projector on divergence free vector fields.\par

The MS system is known to be locally well posed both in the Coulomb gauge and in the Lorentz gauge $\partial_t A_e + \nabla \cdot A = 0$ in sufficiently regular spaces \cite{8r} \cite{9r}, to have weak global solutions in the energy space \cite{7r} and to be globally well posed in a space smaller than the energy space \cite{10r}.\par

A large amount of work has been devoted to the theory of scattering and more precisely to the existence of wave operators for nonlinear equations and systems centering on the Schr\"odinger equation and in particular for the Maxwell-Schr\"odinger system \cite{2r} \cite{4r} \cite{5r} \cite{12r} \cite{14r}. As in the case of the linear Schr\"odinger equation, one must distinguish the short range case from the long range case. In the former case, ordinary wave operators are expected and in a number of cases proved to exist, describing solutions where the Schr\"odinger function behaves asymptotically like a solution of the free Schr\"odinger equation. In the latter case, ordinary wave operators do not exist and have to be replaced by modified wave operators including an additional phase in the asymptotic behaviour of the Schr\"odinger function. In that respect, the MS system in ${I\hskip-1truemm R}^{3+1}$ belongs to the borderline (Coulomb) long range case. We refer to II and \cite{6r} for general background and additional references on that matter.\par

The main step in the construction of the (modified) wave operators consists in solving the local Cauchy problem with infinite initial time. In the long range case where that problem is singular, that step amounts to construct solutions with prescribed (singular) asymptotic behaviour in time. For the MS system in the Coulomb gauge (\ref{1.3e}) (\ref{1.4e}), that step was performed in II by replacing the original system by an auxiliary system, solving the corresponding problem for that system and then returning to the original one. In particular we derived an existence and uniqueness result for solutions of the auxiliary system with prescribed time asymptotics, from which an existence result for solutions of the original system with prescribed time asymptotics follows. However uniqueness was proved only for the auxiliary system, thereby leaving uniqueness for the original one open. The purpose of the present paper is to supplement the previous results with a direct uniqueness result for the original system, expressed in terms of the original functions $(u,A)$.\par

In order to state that result we first replace the equation (\ref{1.4e}) for $A$ by the associated integral equation with prescribed asymptotic data $(A_+, \dot{A}_+)$, namely 
\beq
\label{1.5e}
A = A_0 - \int_t^{\infty} dt' \ \omega^{-1} \sin (\omega (t-t')) P \ {\rm Im} (\overline{u} \nabla_A u) (t')
\eeq

\noi where $\omega = (- \Delta)^{1/2}$ and $A_0$ is the solution of the free wave equation $\sq A_0 = 0$ given by
\beq
\label{1.6e}
A_0 = (\cos \omega t) A_+ + \omega^{-1} (\sin \omega t ) \dot{A}_+ \ .
\eeq

\noi In order to ensure the gauge condition $\nabla \cdot A = 0$, we assume that $\nabla\cdot A_+ = \nabla \cdot \dot{A}_+ = 0$. As a consequence $x \cdot A_0$ is also a solution of the free wave equation. The uniqueness result will be stated for the MS system in the form (\ref{1.3e}) (\ref{1.5e}). Since the Cauchy problem for that system is singular at $t = \infty$, especially as regards the function $u$, the uniqueness result for that system takes a slightly unusual form. Roughly speaking it states that two solutions $(u_i, A_i)$, $i = 1,2$, coincide provided $u_i$ and $A_i - A_0$ do not blow up too fast and provided $u_1 - u_2$ tends to zero in a suitable sense as $t \to \infty$. In particular that result does not make any reference to the asymptotic data for $u$, which should characterize its behaviour at infinity.\par

In order to state the result we need some notation. We denote by $\parallel \, \cdot\, \parallel_r$ the norm in $L^r \equiv L^r ({I\hskip-1truemm R}^3)$, $1 \leq r \leq \infty$ and by $\dot{H}^1 = \dot{H}^1 ({I\hskip-1truemm R}^3)$ the homogeneous Sobolev space
$$\dot{H}^1 = \left \{ v : \nabla v \in L^2\ {\rm and}\ v \in L^6 \right \}\ .$$

\noi We shall need the space
\beq
\label{1.7e}
V_{\star} = \left \{ v :\  <x>^3 v \in L^2\ ,\ <x>^2 \nabla v \in L^2 \right \}
\eeq

\noi where $<\cdot > \ = (1 + |\cdot |^2)^{1/2}$, and the dilation operator
\beq
\label{1.8e}
S = t \partial_t + x \cdot \nabla + \1\ .
\eeq

\noi It follows from the commutation relation $\sq S = (S+2)\sq$ that $SA_0$ satisfies the free wave equation if $A_0$ does. We shall use the notation 
\beq
\label{1.9e}
\widetilde{u}(t) = U(-t) \ u(t)
\eeq

\noi where $U(t) = \exp (i (t/2)\Delta )$ is the unitary group which solves the free Schr\"odinger equation. We denote non negative integers by $j$, $k$, $\ell$. \par

The main result can be stated as follows.\\

\noi {\bf Proposition 1.1.} {\it Let $1 \leq T < \infty$, $I = [T, \infty )$ and $\alpha \geq 0$. Let $A_0$ be a divergence free solution of the free wave equation satisfying 
\beq
\label{1.10e}
\parallel \nabla^k S^j A_0 (t) \parallel_{\infty}\ + \ \parallel \nabla^k x \cdot A_0(t)  \parallel_{\infty}\ \leq C\ t^{-1} \quad \hbox{ for $0 \leq j + k \leq 1$}
\eeq

\noi for all $t \in I$. Let $(u_i, A_i)$, $i = 1,2$, be two solutions of the system (\ref{1.3e}) (\ref{1.5e}) such that $\widetilde{u}_i \in L_{loc}^{\infty} (I, V_{\star})$, $A_i - A_0 \in L_{loc}^{\infty} (I, \dot{H}^1)$ and such that 
\beq
\label{1.11e}
\parallel x^k \nabla^{\ell} \widetilde{u}_i (t) \parallel_2 \ \leq C(1 + \ell n \ t)^{\alpha}\quad  {\it for} \ 0 \leq \ell \leq 1\ , \ 0 \leq k+\ell \leq 3\ ,
\eeq
\beq
\label{1.12e}
\parallel \nabla (A_i - A_0)(t) \parallel_2 \ \leq C \ t^{-1/2} (1 + \ell n \ t)^{2 \alpha}\ ,
\eeq
\beq
\label{1.13e}
\parallel <x/t> (u_1 - u_2)(t)\parallel_2 \ \leq C\ h_{\star}(t) 
\eeq

\noi for all $t \in I$, where $h_{\star} \in {\cal C} (I, {I\hskip-1truemm R}^+)$ is such that the function 
\beq
\label{1.14e}
\overline{h}_{\star}(t) = t(1 + \ell n\ t)^{3+ 9 \alpha} \ h_{\star} (t)
\eeq

\noi be non increasing for $t$ sufficiently large and satisfy 
\beq
\label{1.15e}
\int_t^{\infty} dt'\ t{'}^{-1} \ \overline{h}_{\star}(t') \leq c\ \overline{h}_{\star}(t)
\eeq

\noi for all $t \in I$. \par

Then $(u_1, A_1) = (u_2, A_2)$.}\\

\noi {\bf Remark 1.1.} As mentioned previously, $SA_0$ and $x \cdot A_0$ are solutions of the free wave equation. The time decay in (\ref{1.10e}) is the optimal decay that can be obtained for solutions of that equation. Sufficient conditions on $A_+$, $\dot{A}_+$ ensuring that decay are well known (see for instance \cite{13r}).\\

\noi {\bf Remark 1.2.} Typical functions $h_{\star}$ satisfying the assumptions of the Proposition are
$$h_{\star}(t) = t^{-\lambda} (1 + \ell n\ t)^{\mu}$$

\noi with $\lambda > 1$ and $\mu$ real. \\

\noi {\bf Remark 1.3.} It will be shown below that the solutions of the system (\ref{1.3e}) (\ref{1.5e}) obtained in II (see especially Proposition 7.2 in II) satisfy the assumptions of Proposition 1.1 with $\alpha = 3$ and $h_{\star}(t) = t^{-2} (1 + \ell n\ t)^4$ so that Proposition 1.1 applies to those solutions.  \\

Proposition 1.1 will be proved by going to the above mentioned auxiliary system and generalizing the uniqueness proof for that system obtained in II (see Proposition 4.2 of II).\par

This paper is organized as follows. In Section~2, we derive the auxiliary system which will replace the original system (\ref{1.3e}) (\ref{1.5e}). In Section~3, we collect some notation and preliminary estimates. In Section~4, we derive the uniqueness result, first for the auxiliary system and then for the original one.

\mysection{The auxiliary system}
\hspace*{\parindent}
In this section we perform a change of unknown functions which is well adapted to the study of the system (\ref{1.3e}) (\ref{1.5e}) for large time and we derive the auxiliary system satisfied by the new functions. The unitary group $U(t)$ which solves the free Schr\"odinger equation can be written as 
\beq
\label{2.1e}
U(t) = \exp (i (t/2)\Delta ) = M(t) \ D(t)\ F\ M(t)
\eeq

\noi where $M(t)$ is the operator of multiplication by the function
\beq
\label{2.2e}
M(t) = \exp \left ( ix^2/2t \right )\ ,
\eeq

\noi $F$ is the Fourier transform and $D(t)$ is the dilation operator defined by 
\beq
\label{2.3e}
D(t) = (it)^{-3/2} D_0(t) \qquad , \quad \left ( D_0 (t) f\right ) (x) = f(x/t)\ .
\eeq

\noi We first change $u$ to its pseudo conformal inverse $u_c$ defined by 
\beq
\label{2.4e}
u(t) = M(t)\ D(t) \ \overline{u_c (1/t)}
\eeq

\noi or equivalently
\beq
\label{2.5e}
\widetilde{u}(t) = \overline{F \widetilde{u}_c (1/t)}\ ,
\eeq

\noi where for any function $f$ of space time
$$\widetilde{f}(t, \cdot ) = U(-t) \ f(t, \cdot )\ .$$

\noi Correspondingly we change $A$ to $B$ defined by 
\beq
\label{2.6e}
A(t) = - t^{-1}D_0 (t)\ B(1/t)\ .
\eeq

\noi The transformation $(u, A) \to (u_c, B)$ is involutive. Furthermore it replaces the study of $(u, A)$ in a neighborhood of infinity in time by the study of $(u_c, B)$ in a neighborhood of $t =0$.\par

Substituting (\ref{2.4e}) (\ref{2.6e}) into (\ref{1.3e}) and commuting the Schr\"odinger operator with $MD$, we obtain
\begin{eqnarray*}
&&\left \{ \left (  i \partial_t + (1/2) \Delta_A - g(u) \right  ) u \right \} (t)\\
&&= t^{-2} \ M(t)\ D(t) \left \{ \overline{(i \partial_t + (1/2) \Delta_B - {\check B} - t^{-1} g(u_c))u_c} \right \} (1/t)
\end{eqnarray*}

\noi where for any ${I\hskip-1truemm R}^3$ vector valued function $f$ of space time
\beq
\label{2.7e}
{\check f}(t, x) = t^{-1} x \cdot f(t, x) \ .
\eeq

\noi Furthermore
$${\rm Im} \left ( \overline{u} \nabla_A u\right ) (t) = t^{-3} D_0(t) \left \{ x |u_c|^2 - t^{-1} \ {\rm Im}\ \overline{u}_c \nabla_B u_c \right \} (1/t)$$

\noi by a direct computation, so that the system  (\ref{1.3e}) (\ref{1.5e}) becomes 
$$\hskip 2.1 truecm 
\left \{   \begin{array}{ll}
i \partial_t u_c = - (1/2) \Delta_B u_c + {\check B} u_c + t^{-1} g(u_c) u_c  &\hskip 4 truecm (2.8)\\
& \\
B_2 = {\cal B}_2 (u_c, B) &\hskip 4 truecm (2.9)
 \end{array} \right . 
$$

\noi where $B_0$ is defined by (2.6)$_0$ and
$$B_2 = B - B_0 - B_1 \ , \eqno(2.10)$$
$$B_1 = B_1 (u_c) \equiv - F_1 (P x |u_c|^2)\ , \eqno(2.11)$$
$${\cal B}_2 (u_c, B) \equiv t\ F_2 \left ( P\ {\rm Im}\ \overline{u}_c \nabla_B u_c \right )\ , \eqno(2.12)$$
$$F_j(M) \equiv \int_1^{\infty} d\nu \ \nu^{-2-j}\ \omega^{-1} \sin ( \omega (\nu - 1))D_0(\nu) \ M(t/\nu)\ . \eqno(2.13)$$

\noi Here we take the point of view that $B_1$ is an explicit function of $u_c$ defined by (2.11) and that (2.10) is a change of dynamical variable from $B$ to $B_2$. The equation (2.9) then replaces (\ref{1.5e}).\par

In order to take into account the long range character of the MS system, we parametrize $u_c$ in terms of a complex amplitude $v$ and a real phase $\varphi$ by
$$u_c = v \exp (- i \varphi )\ . \eqno(2.14)$$

\noi The role of the phase is to cancel the long range terms in (2.8), namely the contribution of $B_1$ to ${\check B}$ and the term $t^{-1} g(u_c)$. Because of the limited regularity of $B_1$, it is convenient to split $B_1$ and $B$ into a short range and a long range part. Let $\chi \in {\cal C}^{\infty} ({I\hskip-1truemm R}^3, {I\hskip-1truemm R})$, $0 \leq \chi \leq 1$, $\chi (\xi ) = 1$ for $|\xi | \leq 1$, $\chi (\xi ) = 0$ for $|\xi | \geq 2$. We define
$$\hskip 2.1 truecm 
\left \{   \begin{array}{l}
{\check B}_L = {\check B}_{1L} = F^{\star} \chi (\cdot \ t^{1/2}) F\ {\check B}_1 \\
 \\
 {\check B}_S = {\check B}_0 + {\check B}_{1S} + {\check B}_2 \qquad , \quad {\check B}_{1S} = {\check B}_1 - {\check B}_{1L}\ ,
 \end{array} \right . \eqno(2.15)
$$

\noi We then obtain the following system for $(v, \varphi , B_2)$
$$\hskip 3.5 truecm 
\left \{   \begin{array}{ll}
i \partial_t v = Hv  &\hskip 5.5 truecm (2.16)\\
& \\
\partial_t  \varphi = t^{-1} g(v) + {\check B}_{1L} (v) &\hskip 5.5 truecm (2.17)\\
& \\
B_2 = {\cal B}_2 (v, K) &\hskip 5.5 truecm (2.18)
 \end{array} \right . 
$$

\noi where
$$H \equiv - (1/2) \Delta_K + {\check B}_{S}\ , \eqno(2.19)$$
$$K \equiv B + \nabla \varphi \equiv B + s \ , \eqno(2.20)$$

\noi by imposing (2.17) as the equation for $\varphi$. Under (2.17), the equation (2.8) becomes (2.16). The system (2.16)-(2.18) is the auxiliary system which replaces the original system (1.3) (1.5).

\mysection{Notation and preliminary estimates}
\hspace*{\parindent}
In this section we introduce some notation and collect a number of estimates which will be used throughout this paper. We denote by $\parallel \, \cdot \, \parallel_r$ the norm in $L^r = L^r ({I\hskip-1truemm R}^3)$. For any non negative $k$ we denote by $H^k = H^k ({I\hskip-1truemm R}^3)$ the standard Sobolev spaces 
$$H^k = \left \{ u \in {\cal S}'({I\hskip-1truemm R}^3):\ \parallel u ; H^k \parallel \ = \ \parallel <\omega >^k u\parallel_2\ < \infty \right \}\ ,$$

\noi where $<\,\cdot \, > \ = (1 + |\, \cdot \, |^2)^{1/2}$ and $\omega = (- \Delta )^{1/2}$. In addition we will use the associated homogeneous spaces $\dot{H}^k$ with norm $\parallel u ; \dot{H}^k \parallel \ = \ \parallel \omega^k u\parallel_2$. It will be understood that $\dot{H}^1 \subset L^6$. For any $k \geq 2$ we shall use the notation 
$$\ddot{H}^k = \dot{H}^1 \cap \dot{H}^k\ .$$

\noi For any Banach space $X \subset {\cal S}'({I\hskip-1truemm R}^3)$ we use the notation 
$$FX = \left \{ u \in {\cal S}'({I\hskip-1truemm R}^3): F^{-1} u \in X \right \}\ .$$

\noi For any interval $I$ and any Banach space $X$ we denote by ${\cal C}(I, X)$ the space of strongly continuous functions from $I$ to $X$ and by $L^{\infty} (I, X)$ (resp. $L_{loc}^{\infty} (I, X)$) the space of measurable essentially bounded (resp. locally essentially bounded) functions from $I$ to $X$. For any real numbers $a$ and $b$ we use the notation $a \vee b = \ {\rm Max} (a,b)$ and $a\wedge b =\ {\rm Min} (a,b)$. \par

We next give estimates of the short and long range parts of $B_1$ defined by (2.15), namely 
$$\parallel \omega^m {\check B}_{1S}\parallel\ \leq \ t^{(p-m)/2} \parallel \omega^p {\check B}_{1S} \parallel_2\  \leq \  t^{(p-m)/2} \parallel \omega^p {\check B}_{1} \parallel_2 \eqno(3.1)$$

\noi for $m \leq p$ and similarly
$$\parallel \omega^m {\check B}_{1L}\parallel\ \leq \ \left ( 2t^{-1/2}\right )^{p-m} \parallel \omega^p {\check B}_{1L} \parallel_2\  \leq \  \left ( 2t^{-1/2}\right )^{p-m}  \parallel \omega^p {\check B}_{1} \parallel_2 \eqno(3.2)$$

\noi for $m \geq p$.\par

We now estimate $F_j (M)$ defined by (2.13) (2.3) and $G_j (M)$ defined similarly by~: 
$$G_j (M) = \int_1^{\infty} d \nu \ \nu^{-1-j} \cos (\omega (\nu - 1)) D_0(\nu )\ M(t/\nu)\ . \eqno(3.3)$$

\noi From (2.13) it follows that
$$\omega F_j (M) = F_{j+1}( \omega M)\ , \eqno(3.4)$$
$$\partial_t F_j (M) = F_{j+1} (\partial_t M) \ , \eqno(3.5)$$
$$x \cdot F_j (PM) = F_{j-1} (x \cdot PM)\ . \eqno(3.6)$$

\noi The first two identities are obvious while in (3.6) we have used the identity
$$[x, f(\omega )]\cdot P = 0$$

\noi which holds for any regular function $f$. In addition a direct computation yields
$$x \cdot PM = P (x \otimes M) - 2 \omega^{-2} \nabla \cdot M$$

\noi from which (3.5) can be continued to 
$$x \cdot F_j (PM) = F_{j-1} \left ( P (x \otimes M) - 2 \omega^{-2} \nabla \cdot M\right ) \ . \eqno(3.7)$$

\noi Clearly the identities (3.4) (3.5) (3.6) (3.7) hold with $F_j$ replaced by $G_j$. The following lemma provides an expression for the time derivative of $F_j(M)$ which does not contain the time derivative of $M$.\\

\noi {\bf Lemma 3.1.} {\it Let $F_j(M)$ and $G_j(M)$ be defined by (2.13) and (3.3) respectively. Then}
$$t \partial_t F_j(M) = - F_j \left ( (x \cdot \nabla + j+1)M \right ) + G_j (M)\ . \eqno(3.8)$$
\vskip 5 truemm

\noi {\bf Proof.} From (3.5) we can write 
$$t \partial_t F_j(M) = - \int_1^{\infty}  d\nu\ \nu^{-2-j}\ \omega^{-1}\sin (\omega (\nu - 1)) D_0(\nu)\ \nu \partial_{\nu} M (t/\nu )\ . \eqno(3.9)$$

\noi Using the commutator identity
$$\left ( \nu \partial_{\nu} + x \cdot \nabla \right ) D_0(\nu) = D_0(\nu) \  \nu \partial_{\nu}$$

\noi we obtain
$$t \partial_t F_j(M) = - \int_1^{\infty}  d\nu\ \nu^{-1-j}\ \omega^{-1}\sin (\omega (\nu - 1)) \partial_{\nu}\left ( D_0(\nu)\ M(t/\nu)\right )  - F_j (x \cdot \nabla M)$$

\noi from which (3.8) follows by integration by parts over the $\nu$ variable. \par
\nobreak \hfill $\sq$ \par

In order to estimate $F_j$ and $G_j$ we define
$$I_j (f) (t) = \int_1^{\infty}  d\nu\ \nu^{-j-3/2}\ f(t/\nu ) \eqno(3.10)$$

\noi for any $j \in {I\hskip-1truemm R}$ and for any non negative function $f$ defined in ${I\hskip-1truemm R}^+$. The estimates on $F_j$ and $G_j$ are summarized in the following lemma.\\

\noi {\bf Lemma 3.2.} {\it For any $m, j \in {I\hskip-1truemm R}$ the following estimates hold~:
$$\hbox{\it (1)} \quad \parallel \omega^m F_j(M) \parallel_2 \ \leq \ c \ I_{j+m-2}\left ( \parallel \omega^{m-1} M \parallel_2 \ \wedge\ \parallel \omega^m M \parallel_2\right )\ ,  \  \qquad \eqno(3.11)$$
$$\parallel \omega^m G_j(M) \parallel_2 \ \leq \ c \ I_{j+m-2} \left ( \parallel \omega^{m} M \parallel_2\right ) \ . \eqno(3.12)$$
$$\hbox{\it (2)} \quad \parallel \omega^m x \cdot F_j(PM) \parallel_2 \ \leq \ c \ I_{j+m-3} \left ( \parallel <x> \omega^{m-1} M \parallel_2 \right )   \ , \qquad\qquad \  \eqno(3.13)$$
$$\parallel \omega^m x \cdot G_j(PM) \parallel_2 \ \leq \ c \ I_{j+m-3} \left ( \parallel x\omega^{m} M \parallel_2\ + \ \parallel \omega^{m-1} M \parallel_2\right ) \ . \eqno(3.14)$$

(3) \quad For any $r$, $2 \leq r \leq 4$, 
$$\parallel F_j(M)\parallel_r\ \leq\ c  \int_1^{\infty}  d\nu\ \nu^{-1+2/r} \ \nu^{-j+1/r}  \parallel M(t/ \nu)\parallel_{r_1}   \eqno(3.15)$$

\noi and
$$\parallel G_j(M)\parallel_r\ \leq\ c  \int_1^{\infty}  d\nu\ \nu^{-1+2/r} \ \nu^{-j+1/r}  \parallel \omega M(t/ \nu)\parallel_{r_1} \eqno(3.16)$$

\noi with $3/r_1 = 2 + 1/r$.}\\

\noi {\bf Proof.} \underline{Part (1).} From the definition of $F_j$ and $G_j$, from (3.4) and the analogue for $G_j$, from the identity
$$\parallel \omega^{m} \ D_0(\nu ) v \parallel_2\ = \ \nu^{-m+3/2}  \parallel \omega^{m} v \parallel_2\ ,$$

\noi and from the estimates
$$|\sin ( \omega ( \nu - 1))| \leq 1 \wedge \omega \nu  \qquad , \quad |\cos ( \omega ( \nu - 1))| \leq 1$$

\noi we obtain easily (3.11) and (3.12). \\

\noi \underline{Part (2)} is an immediate consequence of (3.7), of the analogue for $G_j$, and of Part (1). \\

\noi \underline{Part (3).} From the pointwise estimate \cite{1r} \cite{11r}
$$\parallel\sin ( \omega ( \nu - 1))v\parallel_r\ \vee \ \parallel \cos ( \omega ( \nu - 1))v\parallel_r\ \leq  c (\nu - 1)^{-1+4/r} \ \parallel  \omega^{2-4/r}v\parallel_{\overline{r}}$$
\noi with $2 \leq r < \infty$ and $1/r + 1/\overline{r} = 1$, it follows that
$$\parallel F_j(M)\parallel_r\ \leq\ c  \int_1^{\infty}  d\nu( \nu - 1)^{-1+2/r} \ \nu^{-j+1/r}  \parallel \omega^{1-4/r} M(t/ \nu)\parallel_{\overline{r}}$$

\noi and
$$\parallel G_j(M)\parallel_r\ \leq\ c  \int_1^{\infty}  d\nu( \nu - 1)^{-1+2/r} \ \nu^{-j+1/r}  \parallel \omega^{2-4/r} M(t/ \nu)\parallel_{\overline{r}}$$

\noi which imply (3.15) and (3.16) by Sobolev inequalities.\par \nobreak \hfill $\sq$\par

In order to take into account the time decay of norms of some variables as $t$ tends to zero, we shall introduce a function $h \in {\cal C}(I, {I\hskip-1truemm R}^+)$, where $I = (0, \tau ]$ for $0 < \tau \leq 1$, such that the function $\overline{h}(t) \equiv t^{-1} (1 - \ell n\ t)^{\gamma} h(t)$ with $\gamma \geq 0$ be non decreasing in $I$ and satisfy 
$$\int_0^t dt'\ t{'}^{-1} \ \overline{h}(t') \leq c\ \overline{h}(t)$$

\noi for some $c > 0$ and for all $t \in I$. By an elementary computation we then obtain
$$I_j \left ( t^{-\lambda} (1 - \ell n\ t)^\mu \ h\right ) (t) = t^{-1/2-j} \int_0^t dt'\ t{'}^{j-1/2- \lambda }(1 - \ell n\ t')^\mu \ h(t')$$
$$\leq c\ t^{-\lambda} (1 - \ell n\ t)^\mu \ h(t) \eqno(3.17)$$

\noi for any real $\mu$, provided $j + 3/2 > \lambda$. \par

In all the estimates in this paper we denote by $C$ a constant depending on the unknown functions through the available norms. Absolute constants, denoted by $c$ in this section, will in general henceforth be omitted. The letters $j$, $k$, $\ell$ will always denote non negative integers.

\mysection{Uniqueness}
\hspace*{\parindent} 
In this section we prove Proposition 1.1. This will be done by replacing the original system (\ref{1.3e}) (\ref{1.5e}) by the auxiliary system (2.16)-(2.18) and deriving first a uniqueness result for the latter. We recall that the functions $B_1$ and ${\cal B}_2$ are defined (cf. (2.11) (2.12)) by 
$$B_1(v) \equiv - F_1 \left ( P x|v|^2\right )\ , \eqno(4.1)$$
$${\cal B}_2(v, K) = t\ F_2 \left ( P\ {\rm Im} \ \overline{v} \nabla_K v \right ) \ . \eqno(4.2)$$

\noi The latter will be used in general with 
$$K = B + s = B_0 + B_1 (v) + B_2 + \nabla \varphi\ . \eqno(4.3)$$

\noi We shall need the space
$$V = \left \{ v: v\in H^3\ \hbox{and} \ x v \in H^2\right \} \eqno(4.4)$$

\noi with the natural norm, and for $0 < \tau \leq 1$, $I = (0, \tau ]$ and $\alpha \geq 0$, we shall make use of the assumption\\

\noi (A$_+ \alpha)$ \quad  $v \in L_{loc}^{\infty} (I, V)$ and
$$\parallel v(t) ; V \parallel \ \leq C\ L^{\alpha} \eqno(4.5)$$

\noi for all $t \in I$, where $L = 1 - \ell n\ t$.\\

We first prepare the uniqueness result for the system (2.16)-(2.18) with two lemmas.\\

\noi {\bf Lemma 4.1.} {\it Let $0 < \tau \leq 1$, $I = (0, \tau ]$, $\alpha \geq 0$ and let $v$ satisfy $(A_+ \alpha )$. Then\par

(1) $B_1 (v) \in L_{loc}^{\infty} (I, \ddot{H}^4)$, $\nabla {\check B}_1 (v) \in L_{loc}^{\infty} (I, \ddot{H}^2)$, $g \in L_{loc}^{\infty} (I, \ddot{H}^5)$, $\partial_tB_1(v) \in L_{loc}^{\infty} (I, \ddot{H}^2)$ and the following estimates hold for all $t \in I$~:
$$\parallel \nabla^k B_1(v)\parallel_2\ \leq C\  L^{2 \alpha} \qquad \hbox{for $1 \leq k \leq 4$}\ , \eqno(4.6)$$
$$\parallel \nabla^k {\check B}_1(v)\parallel_2\ \leq C\ t^{-1}\ L^{2 \alpha} \qquad \hbox{for $2 \leq k \leq 3$}\ ,  \eqno(4.7)$$
$$\parallel \nabla^k g(v)\parallel_2\ \leq C\  L^{2 \alpha} \qquad \hbox{for $1 \leq k \leq 5$} \ , \eqno(4.8)$$
$$\parallel \nabla^k \partial_t B_1(v)\parallel_2\ \leq C\ t^{-1}\ L^{2 \alpha} \qquad \hbox{for $1 \leq k \leq 2$}\ .  \eqno(4.9)$$

\noi Let in addition $\varphi$ satisfy (2.17). Then $\nabla \partial_t \varphi \in L_{loc}^{\infty} (I, \ddot{H}^2)$ and 
$$\parallel \nabla^{k+1}  \partial_t  \varphi \parallel_2 \ \equiv \ \parallel \nabla^k \partial_t  s\parallel _2\ \leq C\ t^{-1}\ L^{2 \alpha}\quad \hbox{for $1 \leq k \leq 2$} \ .\eqno(4.10)$$

\noi Let in addition $\nabla \varphi (t_0) \in \ddot{H}^2$ for some $t_0 \in I$. Then $\nabla \varphi  \in {\cal C}(I, \ddot{H}^2)$ and 
$$\parallel \nabla^{k+1}  \varphi \parallel_2 \ \equiv \ \parallel \nabla^k s\parallel _2\ \leq C\ L^{1+ 2 \alpha}\quad \hbox{for $1 \leq k \leq 2$} \ .\eqno(4.11)$$

(2) Let in addition $u = v \exp (- i \varphi )$. Then $u$ satisfies $(A_+ \alpha_1)$ with $\alpha_1 = 3 + 7 \alpha$.\par

(3) Let in addition $B_0$ satisfy
$$\parallel \nabla^k B_0(t)\parallel_{\infty}\ \leq C\ t^{-k}\qquad \hbox{for $0 \leq k \leq 1$} \eqno(4.12)$$

\noi and let $B_2 \in L^{\infty} (I, \dot{H}^1)$ satisfy
$$\parallel \nabla B_2(t)\parallel_2 \ \leq C\ L^{2 \alpha} \eqno(4.13)$$

\noi for all $t \in I$. Then ${\cal B}_2 (v, K) \in L_{loc}^{\infty} (I,H^2)$, ${\check {\cal B}}_2(v, K) \in L_{loc}^{\infty} (I, \ddot{H}^2)$, $\partial_t {\cal B}_2 (v, K) \in L_{loc}^{\infty} (I,H^1)$ and the following estimates hold for all $t \in I$~: 
$$\parallel \nabla^k {\cal B}_2(v, K)\parallel_2\ \leq C\  L^{2 \alpha} \qquad \hbox{for $0 \leq k \leq 2$}\ ,  \eqno(4.14)$$
$$\parallel \nabla^k {\check {\cal B}}_2(v, K )\parallel_2\ \leq C\ t^{-1}\ L^{2 \alpha} \qquad \hbox{for $1 \leq k \leq 2$}\ ,  \eqno(4.15)$$
$$\parallel \nabla^k \partial_t {\cal B}_2(v, K)\parallel_2\ \leq C\ t^{-1}\ L^{2 \alpha} \qquad \hbox{for $0 \leq k \leq 1$} \eqno(4.16)$$

\noi where $K$ is given by (4.3).}\\

\noi {\bf Remark 4.1.} The condition $\nabla f \in \ddot{H}^2$ seems to leave some ambiguity on the nature of $f$. However it implies that $\nabla f \in L^{\infty}$ by Sobolev  inequalities and therefore that $<x>^{-1}f\in L^{\infty}$. This occurs in particular in Part (1) for ${\check B}_1(v)$, $\partial_t  \varphi$ and $\varphi$ for fixed time.\\

\noi {\bf Proof.}

\noi \underline{Part (1).} We first derive the estimates (4.6)-(4.11). \par

It follows from (4.1) and (3.11) that 
$$\parallel \nabla^k B_1(v)\parallel_2\ \leq I_{k-1} \left ( \parallel \nabla^{k-1}  x |v|^2\parallel_2 \right ) \ \leq \ C\  L^{2 \alpha} \qquad \hbox{for $1 \leq k \leq 4$}$$

\noi by $(A_+ \alpha)$ and H\"older and Sobolev inequalities. Similarly from (3.13) 
$$\parallel \nabla^k {\check B}_1(v)\parallel_2\ \leq I_{k-2} \left ( \parallel <x> \nabla^{k-1}  x |v|^2\parallel_2 \right ) \ \leq \ C\  L^{2 \alpha} \qquad \hbox{for $2 \leq k \leq 3$}\ . $$

\noi (4.8) is obvious. It follows from (3.8) (3.11) (3.12) that 
\begin{eqnarray*}
\parallel \nabla^k t \partial_t B_1(v)\parallel_2 &\leq& I_{k-1} \left ( \parallel\nabla^{k-1}  (x \cdot \nabla + 2) x |v|^2\parallel_2 + \ \parallel\nabla^{k}  x |v|^2\parallel_2\right ) \\
&\leq& C\  L^{2 \alpha} \qquad \hbox{for $1 \leq k \leq 2$}\ . 
\end{eqnarray*}

\noi (4.10) follows from (2.17) (4.7) (4.8) while (4.11) follows from (4.10) by integration over time.\par

In order to complete the proof, we need to estimate a lower norm of $B_1$, $\nabla {\check B}_1$, $g$, $\partial_t B_1$ and $\nabla \partial_t \varphi$ in order to show that those quantities belong to $\dot{H}^1$. We estimate them in $L^4$ norm by using the special case $r= 4$ of (3.15) (3.16), namely 
$$\parallel F_j(M)\parallel_4\ \leq\   \int_1^{\infty}  d\nu ( \nu-1)^{-1/2} \ \nu^{-j+1/4}  \parallel M(t/ \nu)\parallel_{4/3} \eqno(4.17)$$

\noi and similarly for (3.16), and by using the Hardy-Littlewood-Sobolev (HLS) inequality for $g$. The right hand side of (4.17) and of the other estimates with the appropriate $M$ is then estimated by the use of $(A_+ \alpha )$.\\

\noi \underline{Part (2)} follows from $(A_+ \alpha )$ and (4.11). The required estimates use only the norm of $\nabla \varphi$ in $\ddot{H}^2$ and the worst contribution comes from 
$$\parallel v | \nabla \varphi|^3\parallel_2 \ \leq \ \parallel v \parallel_{\infty}\ \parallel \nabla \varphi\parallel_6^3 \ \leq C\ L^{\alpha + 3(1 + 2 \alpha )}$$

\noi in the estimate of $\parallel \nabla^3 u \parallel_2$.\\

\noi \underline{Part (3).} We first derive the estimates (4.14)-(4.16). We rewrite (4.2) as 
$${\cal B}_2(v, K) = t\ F_2 \left ( P \left ( {\rm Im}\ \overline{v} \nabla v - \left ( B_0 + B_1(v) + B_2 + s \right ) |v|^2 \right ) \right ) \ . \eqno(4.18)$$

\noi From (3.11), we estimate
$$\parallel \nabla^k {\cal B}_2(v, K)\parallel_2\ \leq \ t\ I_k \Big ( \parallel \overline{v} \nabla v \parallel_2\ + \ \parallel B_0 + B_1(v) + s \parallel_{\infty} \ \parallel v \parallel_4^2$$
$$+ \parallel B_2 \parallel_6 \ \parallel v \parallel_6^2 \Big ) \ \leq C\  t\ L^{1+4 \alpha} \qquad \hbox{for $0 \leq k \leq 1$}\ ,  \eqno(4.19)$$
$$\parallel \nabla^2 {\cal B}_2(v, K)\parallel_2\ \leq \ t\ I_2\Big ( \parallel \nabla (\overline{v} \nabla v) \parallel_2\ + \ 2\parallel B_0 + B_1(v) + s \parallel_{\infty} \ \parallel \overline{v} \nabla v \parallel_2$$
$$+ \parallel \nabla (B_0 + B_1(v)) \parallel_{\infty} \ \parallel v \parallel_4^2\ + \ \parallel \nabla (B_2 + s) \parallel_2 \left ( \parallel v  \parallel_{\infty}^2 \ + \ 2\parallel \overline{v} \nabla v \parallel_3\right )  \Big )  \ \leq C\  L^{2 \alpha}   \eqno(4.20)$$

\noi by $(A_+ \alpha )$ (4.6) (4.11) (4.12) (4.13). This proves (4.14). Note that in (4.20) the dominant contribution comes from the term with $\nabla B_0$. All the other terms contribute at most $CtL^{1 + 4 \alpha}$ as in (4.19). The proof of (4.15) is similar, with the factor $t$ omitted, with $I_k$ replaced by $I_{k-1}$ and the factor $x$ absorbed by $v$.\par

(4.16) follows from (4.2) (3.8) (3.11) (3.12). We obtain
$$\parallel \nabla^k \partial_t {\cal B}_2(v, K)\parallel_2\ \leq \ t\ I_k \left  ( \parallel (x\cdot \nabla + 2) {\rm Im}\ \overline{v} \nabla_K v \parallel_2\ + \ \parallel  \nabla {\rm Im}\ \overline{v} \nabla_K v\parallel_2\right ) $$
$$\leq C\  L^{2 \alpha} \qquad \hbox{for $0 \leq k \leq 1$}\ ,  \eqno(4.21)$$

\noi by $(A_+ \alpha )$ (4.6) (4.11) (4.12) (4.13). The dominant contribution comes from 
$$\parallel x \cdot \left ( \nabla B_0\right ) |v|^2 \parallel_2 \ \leq \ \parallel \nabla B_0 \parallel_{\infty}\ \parallel x |v|^2 \parallel_2 \ \leq  C\ t^{-1} \ L^{2 \alpha}\ .$$

\noi In order to complete the proof, in the same way as in Part (1), we estimate the $L^4$ norm of ${\check B}_2(v, K)$ by using (4.17) with the appropriate $M$ and estimating the right hand side thereof through $(A_+ \alpha )$ (4.6) (4.11) (4.12) (4.13). \par \nobreak \hfill $\sq$\par

\noi {\bf Remark 4.2.} For $k = 1$, we have in fact obtained the better estimate 
$$\parallel \nabla {\cal B}_2(v, K)\parallel_2\ \vee\ t \parallel \nabla {\check {\cal B}}_2(v, K)\parallel_2\ \leq C\  t\ L^{1 + 4 \alpha}   \eqno(4.22)$$

\noi in (4.14) (4.15). For $k = 2$, we could also have obtained better estimates by replacing the assumption (4.12) by
$$\parallel \nabla B_0 \parallel_2 \ \leq C\ t^{-1/2}$$

\noi which is also satisfied if $A_0$ is a sufficiently regular solution of the free wave equation. However the estimates (4.14) (4.15) are sufficient for later purposes. \\

We next estimate the difference of two solutions of the auxiliary system (2.16)-(2.18). For two functions or operators of the same nature $f_i$, $i = 1,2$, we shall use the notation $f_\pm = (1/2) (f_1 \pm f_2)$, so that $f_1 = f_+ + f_-$, $f_2 = f_+ - f_-$ and $(fg)_\pm = f_+ g_\pm + f_- g_\mp$. If $(v_i , \varphi_i, B_{2i})$, $i = 1,2$, are two solutions of the auxiliary system (2.16)-(2.18), then $(v_-, \varphi_-, B_{2-})$ satisfies the system
$$i \partial_t v_- = H_+ v_- + H_- v_+ \eqno(4.23)$$
$$\partial_t \varphi_- = t^{-1} g_- + {\check B}_{1L-} \eqno(4.24)$$
$$B_{2-} = t\ F_2 \left ( P \left ( 2 \ {\rm Im} \ \overline{v}_+ \nabla_{K_+} v_- - K_- \left ( |v_+|^2 + |v_-|^2 \right ) \right ) \right . \eqno(4.25)$$

\noi where
$$H_+ = - (1/2) \Delta_{K_+} + (1/2) K_-^2 + {\check B}_{S+}\ , \eqno(4.26)$$
$$H_- = i K_- \cdot  \nabla_{K_+} + (i/2) (\nabla \cdot s_-) + {\check B}_{S-}\ , \eqno(4.27)$$
$$B_{1-} = (1/2) \left ( B_1(v_1) - B_1(v_2)\right ) = - F_1 \left ( 2 P\ {\rm Re} \ x \overline{v}_+  v_-  \right )\ , \eqno(4.28)$$

\noi ${\check B}_{S\pm}$ and ${\check B}_{L\pm}$ are defined by similar formulas, and $g_-$ and $K_\pm$ are obtained from $g_i = g(v_i)$ and
$$K_i = B_0 + B_1 (v_i) + B_{2i} + \nabla \varphi_i \ . \eqno(4.29)$$

For $0 < \tau \leq 1$, $I = (0 , \tau ]$ and $h \in {\cal C} (I, {I\hskip-1truemm R}^+)$, we introduce the assumption\\

\noi (A$_-h)$ \quad $<x> v_- \in L^{\infty} (I, L^2)$ and
$$\parallel <x> v_- (t) \parallel_2 \ \leq C\ h(t) \eqno(4.30)$$

\noi for all $t \in I$.  \\

\noi {\bf Lemma 4.2.} {\it Let $0 < \tau \leq 1$, $I = (0, \tau ]$, $\alpha \geq 0$, and let $h \in {\cal C}(I,{I\hskip-1truemm R}^+ )$ satisfy
$$\int_0^{\tau} dt \ t^{-3/2}\ L^{\alpha}\ h(t) < \infty \ . \eqno(4.31)$$

(1) Let $v_i$, $i = 1,2$ satisfy $(A_+ \alpha )$ with $v_-$ satisfying $(A_-h)$. Then $B_{1-} \in L^{\infty} (I, \dot{H}^1)$, ${\check B}_{1-}  \in L_{loc}^{\infty} (I, \dot{H}^1)$, $g_- \in L_{loc}^{\infty} (I, \ddot{H}^3)$, and the following estimates hold for all $t \in I$~:
$$\parallel \nabla B_{1-}\parallel_2\ \leq C\  I_0 \left ( \parallel v_- \parallel_2 \ L^{\alpha}\right ) \ ,  \eqno(4.32)$$
$$\parallel \nabla {\check B}_{1-}\parallel_2\ \leq C\ t^{-1}\ I_{-1}  \left ( \parallel <x> v_- \parallel_2 \ L^{\alpha}\right ) \ ,\eqno(4.33)$$
$$\parallel \nabla^{k+1}  g_- \parallel_2\ \leq C\left ( \parallel v_- \parallel_2\ +\  \delta_{k,2} \parallel \nabla v_- \parallel_2\right )   L^{ \alpha} \quad \hbox{for $0 \leq k \leq 2$}\ .  \eqno(4.34)$$

\noi Let in addition $\varphi_i$, $i=1,2$, satisfy (2.17) with $v = v_i$. Then $ \partial_t \varphi_- \in L_{loc}^{\infty} (I, \ddot{H}^3)$ and the following estimates hold for all $t \in I$ 
$$\parallel \nabla^{k+1}  \partial_t  \varphi_- \parallel_2 \ = \ \parallel \nabla^k \partial_t  s_-\parallel _2\ \leq \ C\Big \{  \left ( \parallel v_- \parallel_2\ + \ \delta_{k,2} \parallel \nabla v_- \parallel_2\right )  t^{-1}\ L^{\alpha}$$
$$+ \ t^{-1-k/2}\ I_{-1}  \left ( \parallel <x> v_- \parallel_2 \ L^{\alpha}\right ) \Big \} \quad \hbox{for $0 \leq k \leq 2$} \ .\eqno(4.35)$$

(2) Let $B_0$ satisfy 
$$\parallel \nabla^k (t \partial_t)^j B_0\parallel_{\infty}\ \vee \ \parallel \nabla^k {\check B}_0 \parallel_{\infty}\  \leq C\ t^{-k}\qquad \hbox{for $0 \leq j+ k \leq 1$} \eqno(4.36)$$

\noi for all $t \in I$. Let $(v_i, \varphi_i, B_{2i})$, $i =1,2$, be two solutions of the system (2.16)-(2.18) such that $v_i$ satisfy $(A_+ \alpha )$, such that $B_{2i} \in L_{loc}^{\infty} (I, \dot{H}^1)$ with
$$\parallel \nabla B_{2i}(t)\parallel_2 \ \leq C\ L^{2 \alpha} \eqno(4.37)$$

\noi for all $t \in I$, and such that $\nabla \varphi_i (t_0) \in \ddot{H}^2$ with $\nabla \varphi_- (t_0) \in L^2$ for some $t_0 \in I$, so that $s_- = \nabla \varphi_- \in {\cal C}(I,H^2)$ by (4.35). Then the following estimates hold~:}
$$\left | \partial_t \parallel  v_- \parallel_2 \right | \ \leq \ C \Big \{ \parallel \nabla B_- \parallel _2 \ L^{1 + 3 \alpha } \ + \ \parallel \nabla s_- \parallel _2\ L^{\alpha} \ + \ \parallel  s_- \parallel _2\ L^{1 + 3 \alpha}$$
$$+\ \parallel \nabla {\check B}_{1-}\parallel _2\ t^{1/2}\ L^{\alpha}\ + \ \parallel  \nabla {\check B}_{2-} \parallel _2 \ L^{\alpha} \Big \} \ \equiv E(t)\ , \eqno(4.38)$$
$$\left | \partial_t \parallel  x v_- \parallel_2 \right | \ \leq \ \parallel \nabla_{K_+} v_- \parallel _2\ + E(t)\ ,\eqno(4.39)$$
$$\left | \partial_t \parallel  \nabla_{K_+} v_- \parallel_2 \right | \ \leq \ C \Big \{ \left ( \parallel v_- \parallel _2\ + \  \parallel v_- \parallel _3\right ) \ t^{-1} \ L^{2\alpha } \ + \ \parallel s_- \parallel _2\ t^{-1}\ L^{\alpha}$$
$$ + \ \parallel  \nabla B_- \parallel _2\ t^{-1}\ L^{\alpha}\ +\ \parallel  \nabla s_- \parallel _2\ L^{1+ 3\alpha}\ + \parallel  \nabla \nabla \cdot s_- \parallel _2 \ L^{\alpha}$$
$$+ \ \parallel \nabla {\check B}_{1-}\parallel _2\ L^{\alpha}\ + \ \parallel  \nabla {\check B}_{2-} \parallel _2 \ L^{1+ 3\alpha} \Big \} \ , \eqno(4.40)$$
$$ \parallel  \nabla B_{2-} \parallel _2\ \leq t\ I_1 \left ( \parallel v_- \parallel _2\ L^{1 + 3 \alpha} + \left ( \parallel s_- \parallel _2\ + \ \parallel \nabla B_- \parallel _2\right ) L^{2\alpha}\right ) \ , \eqno(4.41)$$
$$ \parallel  \nabla {\check B}_{2-} \parallel _2\ \leq I_0 \left ( \parallel v_- \parallel _2\ L^{1 + 3 \alpha} + \left ( \parallel s_- \parallel _2\ + \ \parallel \nabla B_- \parallel _2\right ) L^{2\alpha}\right ) \ . \eqno(4.42)$$
\vskip 5 truemm

\noi {\bf Remark 4.3.} The assumption that $v_-$ satisfies $(A_- h)$ with $h$ satisfying (4.31) serves to ensure the finiteness of the RHS of (4.33) and is never used otherwise. Similarly the assumption that $\nabla \varphi_- (t_0) \in L^2$ serves only to ensure that $s_- \in {\cal C} (I, L^2)$.\\

\noi {\bf Proof.}

\noi \underline{Part (1).} We first derive the estimates (4.32)-(4.35). It follows from (4.28) (3.11) (3.13) and $(A_+ \alpha )$ that 
$$\parallel  \nabla B_{1-} \parallel _2\ \leq\ 2\ I_0 \left ( \parallel x \overline{v}_+ v_- \parallel _2\right ) \ \leq \ C\ I_0  \left ( \parallel v_- \parallel _2 \  L^{\alpha}\right ) \ ,$$
$$\parallel  \nabla {\check B}_{1-} \parallel _2\ \leq 2\ t^{-1}\ I_{-1} \left ( \parallel <x>^2 \overline{v}_+v_- \parallel _2\right ) \   \leq \ C\ t^{-1}\  I_{-1}  \left ( \parallel <x> v_- \parallel _2\ L^{\alpha}\right ) \ ,$$

\noi while
$$\parallel  \nabla^{k+1} g_- \parallel _2\ = 2 \parallel \nabla^{k-1} \overline{v}_+ v_-\parallel _2$$

\noi from which (4.34) follows by the use of $(A_+ \alpha )$. (4.35) follows from (2.17) (3.2) (4.33) (4.34).\par

In order to complete the proof, we need to estimate a lower norm of $B_{1-}$, ${\check B}_{1-}$ and $g$. As in the proof of Lemma 4.1, part (1), we estimate the $L^4$ norm of those quantities by using (4.17), the HLS inequality and $(A_+ \alpha )$.\\

\noi \underline{Part (2).} We first note that from (2.18) and Lemma 4.1, part (3), especially (4.14)-(4.16), it follows that $B_{2+} \in L_{loc}^{\infty} (I, H^2)$, ${\check B}_{2+}  \in L_{loc}^{\infty} (I, \ddot{H}^2)$, $\partial_t B_{2+} \in L_{loc}^{\infty} (I, H^1)$ and that the following estimate holds for all $t \in I$
$$\parallel  B_{2+} ; H^2\parallel  \vee\  t \parallel  {\check B}_{2+} ; \ddot{H}^2 \parallel  \vee   \parallel  t \partial_t B_{2+} ; H^1\parallel\ \leq C\ L^{2 \alpha}\ . \eqno(4.43)$$

\noi Together with (4.36) and with Lemma 4.1, part (1), especially (4.6) (4.11), this implies that $K_+ \in L_{loc}^{\infty} (I, \ddot{H}^2)$ and that $K_+$ satisfies the estimate
$$\parallel  K_{+} \parallel_{\infty} \ \leq \ C \parallel  K_{+} ; \ddot{H}^2\parallel \ \leq \ C\ L^{1 + 2 \alpha}\ . \eqno(4.44)$$

\noi We next estimate $\parallel  v_- \parallel_2$. From (4.23) (4.27) (3.1), we obtain 
$$\left | \partial_t \parallel  v_- \parallel_2 \right | \ \leq \ \parallel H_- v_+ \parallel _2$$
$$\leq \ C \Big \{ \parallel \nabla B_- \parallel _2\left ( \parallel \nabla v_+ \parallel _3\ + \ \parallel  K_{+} \parallel_{\infty} \ \parallel  v_{+} \parallel_3\right )\ + \ \parallel  \nabla s_- \parallel _2 \left ( \parallel \nabla v_+ \parallel _3\ + \ \parallel  v_{+} \parallel_{\infty}\right )$$
$$+ \ \left ( \parallel  s_- \parallel _2\ \parallel  K_{+} \parallel_{\infty} \ + \ t^{1/2}\ \parallel  \nabla {\check B}_{1-} \parallel_2\right )   \parallel  v_{+} \parallel_{\infty}\ + \ \parallel  \nabla {\check B}_{2-} \parallel_2\ \parallel  v_{+} \parallel_3 \Big \} \eqno(4.45)$$

\noi from which (4.38) follows by the use of $(A_+ \alpha )$ and (4.44).\par

We next estimate $\parallel  xv_{-} \parallel_2$. From (4.23) and the commutation relation
$$[x, H_+] = \nabla_{K_+}$$

\noi we obtain 
$$\left | \partial_t \parallel  xv_- \parallel_2 \right | \ \leq \ \parallel \nabla_{K_+} v_- \parallel _2\ + \ \parallel x H_-  v_+\parallel _2$$

\noi from which (4.39) follows by estimating the last norm in the same way as in (4.45), with the additional factor $x$ everywhere absorbed by $v_+$. We next estimate $\parallel \nabla_{K_+} v_- \parallel _2$. Taking the covariant gradient of (4.23) yields 
$$i \partial_t \nabla_{K_+} v_- = - (1/2) \nabla_{K_+}\Delta_{K_+} v_- + \left ( (1/2) K_-^2 + {\check B}_{S+} \right ) \nabla_{K_+} v_-$$
$$+ \left ( \partial_t K_+ + K_- \nabla_{K_-} + \nabla  {\check B}_{S+}  \right ) v_- + i K_- \cdot \nabla_{K_+}^2 v_+ + i (\nabla K_-)\cdot \nabla_{K_+} v_+$$
$$+ (i/2) (\nabla \cdot s_-) \nabla_{K_+} v_+ + (i/2) (\nabla \nabla \cdot s_-) v_+ +  {\check B}_{S-}  \nabla_{K_+}v_+ + (\nabla  {\check B}_{S-}) v_+ \eqno(4.46)$$

\noi from which we estimate
$$\left |\partial_t \parallel \nabla_{K_+} v_-\parallel_2\right | \ \leq \ \parallel (\partial_t K_+ + \nabla  {\check B}_{S+})v_-\parallel_2 \ + \ \parallel K_-\cdot  \nabla_{K_+}^2 v_+\parallel_2$$
 $$+\  \parallel \nabla K_- \parallel_2 \left ( \parallel \nabla_{K_+} v_+\parallel_{\infty} \ + \  \parallel K_- v_-\parallel_{\infty}\right ) + \ \parallel \nabla \nabla \cdot s_-\parallel_2\ \parallel v_+\parallel_{\infty}$$
$$+ \ \parallel\nabla {\check B}_{1-}\parallel_2 \left ( t^{1/2}  \parallel \nabla_{K_+}  v_+\parallel_{\infty} \ + \ \parallel v_+\parallel_{\infty} \right )$$
$$+ \ \parallel \nabla  {\check B}_{2-}\parallel_2 \left ( \parallel \nabla_{K_+} v_+\parallel_3\ + \ \parallel v_+\parallel_{\infty} \right ) \ . \eqno(4.47)$$
  
\noi We next estimate the first two terms in the right hand side of (4.47). We estimate
$$\parallel (\partial_t K_+ + \nabla {\check B}_{S+}) v_-\parallel_2\ \leq \  \parallel \partial_t (s_+ + B_0 + B_{1+}) + \nabla ({\check B}_0 + {\check B}_{1+} ) \parallel_{\infty} \ \parallel v_- \parallel_2$$
$$+\ \parallel \partial_t  B_{2+} + \nabla {\check B}_{2+} \parallel_6\ \parallel v_-\parallel_3\ \leq \ C \left ( \parallel v_-\parallel_2\ + \ \parallel v_-\parallel_3 \right ) t^{-1}\ L^{2 \alpha} \eqno(4.48)$$

\noi where we have used (4.7) (4.9) (4.10) (4.36) (4.43), and
$$\parallel K_- \nabla_{K_+}^2 v_+ \parallel_2\ \leq \ \parallel s_- \parallel_3 \Big ( \parallel \nabla^2 v_+ \parallel_6 \ + \ \parallel \nabla (s_+ + B_{2+})\parallel_6\ \parallel v_+ \parallel_{\infty} \Big )$$
$$+ \ \parallel s_- \parallel_2 \Big ( \parallel K_+ \parallel_{\infty} \ \parallel \nabla v_+ \parallel_{\infty}\ + \left ( \parallel \nabla (B_0 + B_{1+}) \parallel_{\infty} \ + \ \parallel K_+ \parallel_{\infty}^2 \right ) \parallel v_+ \parallel_{\infty} \Big )$$
$$+ \ \parallel B_- \parallel_6 \Big ( \parallel  \nabla^2 v_+\parallel_3\ +\ \parallel  K_+ \parallel_{\infty}\ \parallel \nabla v_+ \parallel_3\ + \ \parallel \nabla (B_0 + B_{1+}) \parallel_{\infty} \ \parallel v_+ \parallel_3$$
$$+\ \parallel \nabla (s_+ + B_{2+} ) \parallel_6\ \parallel v_+ \parallel_6\ + \ \parallel K_+\parallel_{\infty}^2 \ \parallel v_+ \parallel_3 \Big )$$
$$\leq C \left \{ \parallel s_- \parallel_3 \ L^{1+ 3 \alpha} + \left ( \parallel s_- \parallel_2 \ + \ \parallel \nabla B_- \parallel_2 \right ) t^{-1} \ L^{\alpha}\right \}  \eqno(4.49)$$

\noi where we have used $(A_+ \alpha )$ (4.6) (4.11) (4.36) (4.43) (4.44). \par

Substituting (4.48) (4.49) into (4.47) and estimating the remaining terms of (4.47) by the use of $(A_+ \alpha )$ and (4.44) yields (4.40).\par

We finally estimate $B_{2-}$. From (4.25) and (3.11) (3.13) we obtain
$$\parallel  \nabla B_{2-} \parallel_2\ \leq t\ I_1 \left ( \parallel  v_- \parallel _2\ \parallel  \nabla_{K_+} v_+ \parallel _{\infty} \ + \ \parallel  s_- \parallel _2 \ \parallel  v_+ \parallel _{\infty}^2\ + \ \parallel B_- \parallel _6 \ \parallel v_+ \parallel_6^2 \right ) \ ,$$

$$\parallel  \nabla {\check B}_{2-} \parallel_2\ \leq \ I_0 \Big ( \parallel  <x> v_- \parallel _2\ \parallel  \nabla_{K_+} v_+ \parallel _{\infty} \ + \ \parallel  s_- \parallel _2 \ \parallel  v_+ \parallel _{\infty}\ \parallel  <x> v_+ \parallel _{\infty}$$
$$+ \ \parallel B_- \parallel _6 \ \parallel v_+ \parallel_6\ \parallel  <x> v_+ \parallel _6 \Big )$$

\noi from which (4.41) (4.42) follow by the use of $(A_+ \alpha )$ and (4.44). \par \nobreak \hfill $\sq$ \par

We now state the uniqueness result for the system (2.16)-(2.18).\\

\noi {\bf Proposition 4.1.} {\it Let $0 < \tau \leq 1$, let $I = (0, \tau ]$, $\alpha \geq 0$ and let $h \in {\cal C} (I, {I\hskip-1truemm R}^+)$ be such that $\overline{h} (t) = t^{-1} (1 - \ell n\ t)^{\alpha} h(t)$ be non decreasing and satisfy  
$$\int_0^t dt'\ t{'}^{-1} \overline{h} (t') \leq c\ \overline{h}(t) \eqno(4.50)$$

\noi for some $c > 0$ and for all $t \in I$. Let $B_0$ satisfy (4.36) for all $t \in I$. Let $(v_i, \varphi_i, B_{2i})$, $i = 1,2$, be two solutions of the system (2.16)-(2.18) such that $v_i$ satisfy $(A_+ \alpha )$, such that $B_{2i}  \in L_{loc}^{\infty} (I, \dot{H}^1)$ and satisfy (4.37) for all $t\in I$, and such that $\nabla \varphi_i (t_0) \in \ddot{H}^2$ for some $t_0 \in I$.  Assume in addition that $\varphi_- (0) = 0$ and that $v_-$ satisfy $(A_- h)$.\par

Then $(v_1, \varphi_1, B_{21}) = (v_2, \varphi_2, B_{22})$.}\\

\noi {\bf Proof.} Note first that (4.50) implies (4.31) so that Lemma 4.2 can be applied. From (4.35) with $k = 0$ and mild assumptions on $v_-$, it follows that $\varphi_-(t)$ has a limit in $\dot{H}^1$ as $t \to 0$, thereby giving a meaning to the assumption $\varphi_- (0) = 0$. Actually it follows from (4.35) (4.50) and $(A_- h)$ that the limit exists in $\ddot{H}^3$. \par

We first prove the proposition for $\tau$ sufficiently small by using Lemma 4.2. We define
$$y_0 = \ \parallel <x> v_- \parallel_2 \qquad , \quad y_1 = \ \parallel \nabla_{K_+} v_- \parallel_2\ ,$$
$$Y_0 = \ \mathrel{\mathop {\rm Sup}_{t \in I}}\ h (t)^{-1} y_0 (t)\ .$$

\noi From Lemma 4.2, especially  (4.32)  (4.33)  (4.35) (4.50) and from  (3.17), we obtain
$$\parallel \nabla B_{1-}\parallel_2\  \leq C\ Y_0\ L^{\alpha} h\ , \eqno(4.51)$$
$$\parallel \nabla {\check B}_{1-}\parallel_2\  \leq C\ Y_0\ t^{-1} \ L^{\alpha} h\ , \eqno(4.52)$$
$$\parallel \nabla^{k+1}  \varphi_- \parallel_2\  = \ \parallel \nabla^k  s_- \parallel_2\ \leq C\ \left \{ Y_0\ t^{-k/2 }  \ L^{\alpha} h + \ \delta_{k,2} \int_0^t dt'\ t{'}^{-1} \ L{'}^{\alpha}\ y_1(t')\right \} \eqno(4.53)$$

\noi for $0 \leq k \leq 2$ and for all $t\in I$, with $L' = 1 - \ell n\ t'$. The time integral of the last term in (4.53) converges because of the estimate
$$\parallel \nabla  v_- \parallel_2^2\  \leq\ \parallel v_- \parallel_2\ \parallel \Delta  v_- \parallel_2 \ \leq C \ Y_0\ L^{\alpha} h$$

\noi and we have replaced the ordinary derivative by the covariant one in that integral, thereby producing an innocuous term with $Y_0L^{1 + 3\alpha} h $. On the other hand from (4.41)  (4.51)  (4.53) (3.17) we obtain
$$
\parallel \nabla  B_{2-} \parallel_2\  \leq \ C\left\{  Y_0 \ t \ L^{1+3\alpha} h + t\ I_1 \left ( \parallel \nabla B_{2-} \parallel_2 L^{2\alpha}\right ) \right \} \ . \eqno(4.54)$$

\noi From the assumptions on $B_{2i}$, it follows that $B_{2-} \in L_{loc}^{\infty} (I, \dot{H}^1)$ with
$$\parallel \nabla  B_{2-}(t) \parallel_2\  \leq \ C\ L^{2\alpha} \ .$$ 

\noi Using that fact, one obtains easily from (4.54) that
$$\parallel \nabla  B_{2-} \parallel_2\  \leq \ C\ Y_0 \ t \ L^{1+3\alpha} \ h \eqno(4.55)$$ 

\noi for all $t \in I$ and for $\tau$ sufficiently small. Substituting that result into  (4.42) yields
$$\parallel \nabla  {\check B}_{2-} \parallel_2\  \leq \ C\ Y_0\ L^{1+3\alpha} \  h \ . \eqno(4.56)$$

\noi Substituting (4.51) (4.52) (4.53) (4.55) (4.56) into (4.40) yields
$$\left | \partial_t y_1\right | \leq C \left \{ Y_0 \ t^{-1} \ L^{2\alpha}\ h + (Y_0\ h \ y_1)^{1/2} \ t^{-1}  \ L^{2\alpha}+ \ L^{\alpha}  \int_0^t dt'\ t{'}^{-1} \ L{'}^{\alpha}\ y_1(t') \right \} \eqno(4.57)$$

\noi which takes the form 
$$\left | \partial_t y\right | \leq f + g\ y^{1/2} + C  \int_0^t dt'\ t{'}^{-1} \ L{'}^{2\alpha}\ y(t') \eqno(4.58)$$

\noi with
$$y_1 = Y_0 \ y\quad , \quad f = C \  t^{-1} \ L^{2\alpha}\  h \quad , \quad g = C\ t^{-1} \\ L^{2\alpha\ }h^{1/2}\ .$$

\noi We define
$$z(t) = \int_0^t dt'\ t{'}^{-1} \ L{'}^{2\alpha}\ y(t')$$

\noi so that $t \partial_t z = L^{2\alpha}y$ and
$$F(t) = \int_0^t dt' \ f(t')\ .$$

\noi Integrating (4.58) over time with $y(0) = 0$ yields 
$$y(t) \leq F(t) +  \int_0^t dt'\ g(t')\ y(t')^{1/2} + C \int_0^t dt' \ z(t')$$
$$\leq F(t) + C \ t\ z(t) + z(t)^{1/2} \left \{ \int_0^t dt'\ t' \ L{'}^{-2 \alpha}\ g^2 (t') \right \}^{1/2}\ . \eqno(4.59)$$

\noi The last integral is estimated by
$$\int_0^t dt'\ t'\ L{'}^{-2 \alpha}\ g^2(t') \leq C\ L^{2\alpha} \ h$$

\noi by (4.50), so that (4.59) yields
$$y(t) \leq F(t) + C \left ( t\ z(t) + (z(t)h(t))^{1/2} \ L^{\alpha} \right ) \eqno(4.60)$$

\noi and therefore
$$\partial_t z \leq t^{-1} \ L^{2 \alpha} \ F(t) + C\ L^{2\alpha}\ z + C(z\ h)^{1/2} \ t^{-1}\ L^{3\alpha} \ . \eqno(4.61)$$

\noi Integrating (4.61) over time (see for instance Lemma 2.3 in \cite{3r}) we obtain
$$z(t) \leq \exp \left ( C\ t\ L^{2 \alpha} \right ) \left \{ \int_0^t dt'\ t{'}^{-1}\ L{'}^{2\alpha} \ F(t') + \left (  \int_0^t dt'\ t{'}^{-1}\ L{'}^{3\alpha} \ h(t')^{1/2} \right )^2 \right \} \ .\eqno(4.62)$$

\noi We next estimate
\begin{eqnarray*}
\int_0^t dt'\ t{'}^{-1}\ L{'}^{2\alpha} \ F(t') &\leq& C\ L^{4\alpha}\ h \ ,\\
\left ( \int_0^t dt'\ t{'}^{-1}\ L{'}^{3\alpha} \ h(t')^{1/2}\right )^2 &\leq& \left ( \int_0^t dt'\ t{'}^{-2}\ L{'}^{\alpha} \ h(t')\right ) \left ( \int_0^t dt'\ L{'}^{5\alpha} \right ) \\
&\leq& C\ L^{6\alpha}\ h 
\end{eqnarray*}

\noi by (4.50) and an elementary computation. Substituting those estimates into (4.62) yields
$$z(t) \leq C\ L^{6\alpha}\ h$$

\noi and therefore by (4.60)
$$y_1 (t) \leq  C\ Y_0 \ L^{4\alpha}\ h \ . \eqno(4.63)$$

\noi Substituting (4.51) (4.52) (4.53) (4.55) (4.56) (4.63)  into (4.38) (4.39) yields
$$\left | \partial_t y_0 \right | \leq C\ Y_0 \ t^{-1/2} \  h$$

\noi and therefore by integration over time with $y_0 (0) = 0$
$$y_0(t)  \leq C\ Y_0 \ t^{1/2}\ h$$

\noi so that
$$Y_0 \leq C\ Y_0\ \tau^{1/2}$$

\noi which implies $Y_0 = 0$ and therefore $v_- = 0$ for $\tau$ sufficiently small. Substituting that result into (4.35) (4.55) shows that $\varphi_- = 0$, $B_{2-} =0$ and therefore $(v_1, \varphi_1, B_{21}) = (v_2, \varphi_2, B_{22})$.\par

The extension of the proof to the case of general $\tau$ proceeds by similar but more standard arguments.\par\nobreak \hfill $\sq$ \par

We now turn to the proof of Proposition 1.1. \\

\noi {\bf Proof of Proposition 1.1.} The first step consists in rewriting that proposition in an equivalent form in terms of the variables $(u_c, B_2)$ where $u_c$ is defined by (2.4) or (2.5). \\

\noi {\bf Proposition 4.2.} {\it Let $0 < \tau \leq 1$, let $I = (0, \tau ]$, $\alpha \geq 0$. Let $A_0$ be a divergence free solution of the free wave equation such that $B_0$ defined by (2.6)$_0$ satisfy (4.36) for all $t\in I$. Let $(u_{ci}, B_{2i})$, $i = 1,2$, be two solutions of the system (2.8) (2.9) such that $u_{ci}$ satisfy $(A_+ \alpha )$, that $(B - B_0)_i \equiv B_1 (u_{ci}) + B_{2i} \in L_{loc}^{\infty} (I, \dot{H}^1)$ with 
$$\parallel \nabla (B -  B_0)_i(t) \parallel_2\  \leq C(1 - \ell n\ t)^{2 \alpha} \eqno(4.64)$$

\noi for all $t \in I$. Assume in addition that $u_{c-}$ satisfy $(A_- h)$ for some function $h \in {\cal C}(I, {I\hskip-1truemm R}^+)$ such that the function
$$\overline{h}(t) = t^{-1} (1 - \ell n\ t)^{3 + 9 \alpha} \ h(t) \eqno(4.65)$$

\noi be non decreasing for $t$ sufficiently small and satisfy
$$\int_0^t dt'\ t{'}^{-1} \ \overline{h}(t') \leq c\ \overline{h}(t) \eqno(4.66)$$

\noi for all $t\in I$. \par

Then $(u_{c1},B_{21}) = (u_{c2},B_{22})$.}\\

We first show the equivalence of Proposition 1.1 and Proposition 4.2, with $\tau = T^{-1}$ and $h(t) = h_{\star} (1/t)$. The equivalence of (1.10) for $A_0$ and (4.36) for $B_0$ follows from (2.6)$_0$ and from the relations
\begin{eqnarray*}
&&SA_0(t) = t^{-1}\ D_0(t)\left ( t \partial_t B_0\right ) (1/t)\\
&&(x\cdot A_0)(t) = - t^{-1}\ D_0(t)\ {\check B}_0 (1/t)\ .
\end{eqnarray*}

The equivalence of the assumption on $u_i$ in Proposition 1.1 with $(A_+ \alpha )$ for $u_{ci}$ follows from (2.5), from the fact that $V = FV_{\star}$ and from the commutation relation
$$x\ U(-t) = U(-t) (x + it \nabla )$$

\noi which implies that 
$$\left | \ \parallel U(-t) v ; V \parallel \ - \ \parallel v;V \parallel\ \right | \leq |t| \parallel v; H^3\parallel$$

\noi so that $(A_+ \alpha )$ for $u_{ci}$ is equivalent to $(A_+ \alpha )$ for $\widetilde{u}_{ci}$. The equivalence of (1.12) for $A-A_0$ with (4.64) for $(B - B_0)$ follows from (2.6). Finally the equivalence of the assumption (1.13) for $u_-$ with $(A_- h)$ for $u_{c-}$ follows from (2.4).\par

We are now left with the task of deriving Proposition 4.2 from Proposition 4.1.\\

\noi {\bf Proof of Proposition 4.2.} \par

Let $(u_{ci}, B_{2i})$ satisfy the assumptions of Proposition 4.2. We need to construct $(v_i, \varphi_i, B_{2i})$ satisfying the assumptions of Proposition 4.1. The main step is to construct $\varphi_i$ from $u_{ci}$. Now it follows from (2.14) that (2.17) is equivalent to 
$$\partial_t \varphi = t^{-1}\ g(u_c) + {\check B}_{1L} (u_c) \eqno(4.67)$$

\noi and we define the phases $\varphi_i$ by integrating that equation over time with some initial condition $\nabla \varphi_i(t_0) \in \ddot{H}^2$. It follows from Lemma~4.1, part (1) with $v$ replaced by $u_c$ that $\varphi_i$ satisfy (4.11) and from Part (2) with $(v, u, \varphi )$ replaced by $(u_c, v, - \varphi )$ that $v_i$ defined by (2.14) satisfy $(A_+ \alpha_1 )$ with $\alpha_1 = 3 + 7 \alpha$. Furthermore, again by Lemma 4.1, part (1), $B_1(u_{ci}) = B_1 (v_i)$ satisfy the assumptions made on $B - B_0$ in Proposition 4.2, so that the latter are equivalent to the assumptions made on $B_2$ in Proposition 4.1. \par

It remains to ensure the assumption $(A_- h_1)$ on $v_-$ for a suitable $h_1$. This will be done by suitably adjusting the initial conditions for the phases $\varphi_i$. If $\varphi_i$, $i = 1,2$, satisfy (4.67), then $\varphi_-$ is estimated in the same way as in Lemma 4.2 with $v$ replaced by $u_c$, namely  
$$\parallel \nabla^{k+1}  \partial_t  \varphi_- \parallel_2 \ \leq C \Big \{\left ( \parallel  u_{c-} \parallel _2 \ + \ \delta_{k,2}  \ \parallel \nabla u_{c-}  \parallel _2 \right ) t^{-1} \ L^{\alpha}$$
$$+ \ t^{-1-k/2}\ I_{-1} \left ( \parallel <x> u_{c-} \parallel _2 \ L^{\alpha} \right ) \Big \}$$
 $$\leq \ C\left  \{  t^{-1}\ L^{\alpha} \left ( h + \delta_{k,2}\ h^{2/3}\ L^{\alpha /3}\right ) + \  t^{-1-k/2}\ L^{\alpha} h \right  ) \eqno(4.68)$$
 
 \noi by using $(A_- \alpha )$ for $u_{c-}$ and estimating 
 $$\parallel \nabla u_{c-} \parallel_2\ \leq\  \parallel u_{c-} \parallel_2^{2/3}\  \parallel \nabla^3 u_{c-} \parallel^{1/3} \ \leq C\ h^{2/3}\ L^{\alpha /3}\ . $$
 
 We now adjust the phases as follows. We choose arbitrarily $\varphi_1$ with $\nabla \varphi_1(t_0) \in \ddot{H}^2$, we define $\varphi_-$ by integrating (4.68) over time with initial condition $\varphi_-(0) = 0$, and we define $\varphi_2 = \varphi_1 - 2 \varphi_-$. The integral of (4.68) converges for $0 \leq k \leq 2$ by (4.66) and yields
 $$\parallel \nabla^{k+1}  \varphi_- \parallel_2 \ \leq C \left  \{\left ( 1 + t^{-k/2}\right ) L^{\alpha} \ h \ + \ \delta_{k,2}  \ \ L^{\alpha /3}\ h^{2/3} \right \} \ .   \eqno(4.69)$$
 
\noi In particular $\nabla \varphi_2 (t_0) \in \ddot{H}^2$ and
$$\parallel \nabla \varphi_- \parallel_2 \ \leq C L^{\alpha} \ h \ . \eqno(4.70)$$

\noi We can now estimate $v_-$.
\begin{eqnarray*}
v_- &=& (1/2) \left ( u_{c1} \exp \left ( i \varphi_1\right ) - u_{c2}\exp \left ( i \varphi_2\right ) \right )\\
&=& (1/2) \exp \left ( i \varphi_1\right )\left ( u_{c1} - u_{c2} + u_{c2} \left ( 1 - \exp \left ( -2i  \varphi_- \right  )\right  )\right .
\end{eqnarray*}

\noi so that 
$$\parallel <x> v_- \parallel_2  \ \leq \ \parallel <x> u_{c-}\parallel_2\ + \ c \parallel u_{c2} \parallel_3 \ \parallel \nabla \varphi_- \parallel_2$$
$$\leq C\ h\left (1 + L^{2\alpha} \right ) \eqno(4.71)$$

\noi by $(A_+ \alpha )$, $(A_- h)$ and (4.70). Therefore $v_-$ satisfies $(A_- h_1)$ with $h_1 = h L^{2 \alpha}$. Finally the assumption (4.50) for
$$\overline{h}_1(t) = t^{-1} (1 - \ell n\ t)^{\alpha_1} \ h_1(t) = t^{-1} (1 - \ell n\ t)^{3 + 9 \alpha} \ h(t)$$
\noi is equivalent to (4.66).\par

Proposition 4.2 then follows from Proposition 4.1. \par\nobreak \hfill $\sq \ \sq$ \par

We conclude this section by showing that Proposition 1.1 applies to the solutions of the system (1.3) (1.5) constructed in II. Because of the equivalence of Proposition 1.1 and Proposition 4.2, it is sufficient to consider the solutions $(u_c, B_2)$ of the system (2.8) (2.9). In II we proved the existence of solutions of that system with prescribed asymptotic behaviour $(u_{ca}, B_{2a})$ as $t$ tends to zero, with $u_c = v \exp (- i \varphi )$ and $u_{ca} = v_a \exp (- i \varphi_a)$. It follows from Lemma 6.1 in II (note that $v$, $v_a$ of this paper are $w$, $w_a$ of II) that $v_a \in L^{\infty} (I, V)$, $\nabla \varphi_a \in L_{loc}^{\infty} (I, \ddot{H}^2)$ and $B_{2a} \in L^{\infty} (I,\dot{H}^1)$ and from Proposition 7.1 in II that the same properties hold for $v$, $\nabla \varphi$ and $B_2$. In particular $v$ satisfies $(A_+ 0)$ so that by Lemma 4.1, part (2), $u_c$ satisfies $(A_+3)$. Furthermore $v - v_a$ satisfies $(A_- h)$ with 
$$h = t^{2} (1 - \ell n\ t)^4 \eqno(4.72)$$

\noi and $(\varphi - \varphi_a) (0) = 0$. In particular if $(v_i, \varphi_i, B_{2i})$, $i=1,2$, are two solutions of the system (2.16)-(2.18) associated with the same $v_a$, then $v_-$ satisfies $(A_- h)$ with $h$ given by (4.72), so that by the same argument as in the proof of Proposition 4.2, but now with $u_c$ and $v$ interchanged, $u_{c-}$ also satisfies $(A_- h)$. Therefore the solutions constructed in II for fixed $(u_{ca}, B_{2a})$ satisfy the assumptions of Proposition 4.2 with $\alpha = 3$ and $h$ given by (4.72), which proves uniqueness for fixed $(u_{ca}, B_{2a})$ as chosen in II. Note that the proof uses Proposition 4.1 with $\alpha_1 = 24$ and
$$h_1 = t^2 (1 - \ell n\ t)^6\ .$$
\vskip 5 truemm

\noi {\bf Remark 4.4.} The construction of $v$ from $u_c$ or of $u_c$ from $v$ through (2.14), (2.17) or (4.67) and Lemma 4.1, part (2) entails some loss in the estimates, typically from $(A_+ \alpha )$ to $(A_+ \alpha_1 )$ with $\alpha_1 = 3 + 7 \alpha$. In the existence proof given in II, starting from $v$ satisfying $(A_+0)$, we obtained $u_c$ satisfying only $(A_+ 3)$, and a direct uniqueness proof for $u_c$ can start only from that weaker assumption. The reconstruction of $v$ with only $u_c$ available to start with then produces another loss and yields only $(A_+24)$ for $v$. Only that weaker assumption can then be used in the auxiliary uniqueness result of Proposition 4.1, even though $v$ was known to satisfy the better estimate $(A_+0)$ to start with.\\

\noi {\bf Remark 4.5.} The loss on $\alpha$ in Lemma 4.1, part (2) can be reduced by using a more accurate estimate than (4.5) in $(A_+ \alpha )$. If one assumes instead 
$$\parallel \nabla^k  <x>^{\ell} v \parallel_2 \ \leq C L^{k\alpha} \qquad \hbox{for $0 \leq \ell \leq 1$, $0 \leq k+ \ell \leq 3\ ,$}$$

\noi then Lemma 4.1, part (2) still holds with $\alpha_1 = 1 + 3 \alpha$.

\end{document}